\input amstex

\def\b1{\text{\bf 1}}

\def\CA{{\Cal A}}

\def\CC{{\Cal C}}
\def\CD{{\Cal D}}
\def\CF{{\Cal F}}

\def\CH{{\Cal H}}

\def\CL{{\Cal L}}
\def\CM{{\Cal M}}

\def\CR{{\Cal R}}

\def\CV{{\Cal V}}

\def\#{\,\check{}}

\def\Ker{\text{Ker}}


\def\hra{\hookrightarrow}
\def\iso{\buildrel\sim\over\rightarrow} 

\def\lra{\longrightarrow}


\documentstyle{conm-p}
\NoBlackBoxes

\topmatter
\title     A remark on primitive cycles and Fourier-Radon transform       \endtitle
\author A.~Beilinson \endauthor
\leftheadtext{A.~BEILINSON}

\address Department of Mathematics, University of Chicago, Chicago, IL 60637\endaddress

\email sasha\@math.uchicago.edu\endemail

\thanks The author was supported in part by NSF
Grant DMS-0401164.\endthanks




\endtopmatter

\document

The aim of this note is to point out that Brylinski's  Radon transform \cite{B} is a natural instrument
for the Green-Griffiths approach to Hodge conjecture \cite{GG}, \cite{BFNP}. In particular, some principal results of  \cite{BFNP} follow from the general fact that Radon transform preserves primitive cohomology (while reversing its grading). As was noticed by Drinfeld, this
assertion is immediate from the basic  Fourier transform functoriality  \cite{L}.\footnote{My initial argument was less elegant (it used relative Lefschetz decomposition).}

This note originates from a talk  given at a student Hodge theory seminar.   I am grateful to V.~Drinfeld for his enlightening  comment, to  D.~Kazhdan for a discussion, and to M.~Kerr and G.~Pearlstein for an exchange of letters. 

\subhead 1. A reformulation of the Hodge conjecture\endsubhead For a compact complex algebraic variety
$X$ let $N_i H^\cdot (X,\Bbb Q)$ be the niveau filtration on its cohomology  (it is Poincar\'e dual to more commonly used coniveau filtration; conjecturally, the two filtrations are complementary). Thus $N_1 H^\cdot (X,\Bbb Q)$ is the intersection of 
kernels of all  restriction maps $H^\cdot (X,\Bbb Q)\to H^\cdot (Y,\Bbb Q)$, where $Y\neq X$ is a closed algebraic subvariety of $X$. According to Totaro and Thomas, see \cite{BFNP} th.~6.5, the Hodge conjecture amounts to the next assertion: For every projective smooth $X$ of dimension $2n$ the subspace of Hodge $(n,n)$-classes in $H^{2n}(X,\Bbb Q)$ has zero intersection with $N_1 H^{2n}(X,\Bbb Q)$. Of course, it suffices to consider the subspace of primitive Hodge classes. Thus every description of $N_1 H^{2n}(X,\Bbb Q)^{\text{prim}}$ provides a reformulation of the Hodge conjecture. The articles \cite{GG} and \cite{BFNP} provide one
such  description; we present it in the last line of the note.

\remark{Remark} As was pointed out by the referee, Kerr and Pearlstein can treat similarly  Grothendieck's generalized Hodge conjecture.\endremark

\remark{Question} For $\gamma$ in a given term of coniveau filtration, what can one say about simplest possible singularities of $Y$ with $\gamma|_Y \neq 0$?  (E.g., by Thomas, for algebraic $\gamma$, i.e., for $\gamma$ in the deepest term of coniveau filtration, the singularities are ODP.)\endremark

\subhead 2. Radon transform {\rm (\cite{B})}\endsubhead We  play with complex algebraic varieties and   $\Bbb Q$-sheaves. An arbitrary ground field and  $\Bbb Q_\ell$-sheaves will do as well.

 For an algebraic variety $Z$, we denote by $D(Z)$ the derived category of bounded constructible $\Bbb Q$-complexes on $Z$; let $\CM (Z)\subset D(Z)$ be  the category of perverse sheaves on $Z$,
 ${}^p\! H : D(Z)\to\CM (Z)$ the cohomology functor (\cite{BBD}). For smooth $Z$ let $\CM^{\text{sm}}(Z)\subset \CM(Z)$ be the Serre subcategory of smooth perverse sheaves (i.e., local systems); it generates the thick subcategory $D^{\text{sm}}(Z)\subset D(Z)$ of complexes with smooth  cohomology. The Verdier quotient $\bar{D}(Z):= D(Z)/D^{\text{sm}}(Z)$ is a t-category with heart $\bar{\CM}(Z):= \CM (Z)/\CM^{\text{sm}}(Z)$. The latter is an Artinian $\Bbb Q$-category; the projection $\CM (Z)\to \bar{\CM}(Z)$ identifies the subcategory of non-smooth irreducible perverse sheaves on $Z$ with that of  irreducible objects in $\bar{\CM}(Z)$.

Let $V$ be a vector space of dimension $n\ge 2$, $V^\vee$ its dual. Let $\Bbb P$, $\Bbb P^\vee$ be the corresponding projective spaces, $i: T\hra \Bbb P \times \Bbb P^\vee $ be the incidence correspondence. Let $p$, $p^\vee $ be the projections $\Bbb P \times \Bbb P^\vee
\rightrightarrows \Bbb P,\,\Bbb P^\vee \!$, and   $p_{(T)}$, $p_{(T)}^\vee$ be  their restrictions to $T$.
The Radon transform functor $\CR : D(\Bbb P)\to D(\Bbb P^\vee )$ is $\CR (M):=  p_{(T)!}^\vee  p_{(T)}^* M [n-2]$. Interchanging $\Bbb P$ and $\Bbb P^\vee$, we get 
$\CR^\vee : D(\Bbb P^\vee )\to D(\Bbb P )$, etc.
Notice that $\CR$ sends $D^{\text{sm}}(\Bbb P)$ to $ D^{\text{sm}}(\Bbb P^\vee )$, so we have $\bar{\CR}: \bar{D}(\Bbb P)\to \bar{D}(\Bbb P^\vee )$.

\proclaim{Theorem {\rm (\cite{B} 3.1)}}
 The compositions $\bar{\CR}\bar{\CR}^\vee$, $\bar{\CR}^\vee \bar{\CR}$ are  Tate twist functors $M\mapsto M(2-n)$. The functors $\bar{\CR}$, $\bar{\CR}^\vee$ are  t-exact, hence they yield equivalences of the abelian categories  $\bar{\CM}(\Bbb P)\leftrightarrows\bar{\CM}(\Bbb P^\vee )$. \qed\endproclaim

\subhead 3. Fourier transform {\rm (\cite{B}, \cite{L})} \endsubhead The formalism of constructible sheaves
extends to algebraic stacks of finite type (\cite{LMB}, \cite{LO}).
The group $\Bbb G_m$ acts on any vector space by homotheties. Consider the quotient stacks $\CV := V/\Bbb G_m$, $\CV^\vee :=V^\vee /\Bbb G_m$, $\CA^1 := \Bbb A^1 /\Bbb G_m$. The open embedding $j_V : V^\circ :=V\smallsetminus \{0\} \hra V$ yields one $j_{\CV} : \Bbb P \hra \CV$, etc. The canonical pairing map $\mu : V\times V^\vee \to \Bbb A^1$  yields  $\mu: \CV \times\CV^\vee \to \CA^1$.
 Let $pr, pr^\vee : \CV\times\CV^\vee \rightrightarrows \CV,\,\CV^\vee $ be the projections.  One has  the (homogenous) Fourier transform  $\CF : D(\CV )\to D(\CV^\vee )$, $\CF (N):= pr^\vee_! (pr^* N\otimes \mu^* j_{\CA^1  *}\Bbb Q  )[n-1]$, see \cite{L} 1.5, 1.9. Interchanging $\CV$ and $\CV^\vee$, we get 
$\CF^\vee : D(\CV^\vee )\to D(\CV )$.

\proclaim{Theorem {\rm (\cite{L} 3.1, 4.2)}}
  The compositions $\CF\CF^\vee$, $\CF^\vee \CF$ are  Tate twist functors $N\mapsto N(-n)$. The functors $\CF$, $\CF^\vee$ are t-exact, hence they yield equivalences of the abelian categories $\CM (\CV )\leftrightarrows \CM (\CV^\vee )$.
 \qed\endproclaim

Consider the closed embeddings $i_V : \{ 0\} \hra V$, $i_{\CV} : B_{\Bbb G_m}= \{0\}/\Bbb G_m \hra \CV$, etc. The projection $j_{\CA^1 *} \Bbb Q \to i_{\CA^1!} \Bbb Q (-1)[-1]$ yields 
a natural morphism $$j_{\CV^\vee}^* \CF j_{\CV !} \to \CR (-1), \tag 1$$ which becomes an isomorphism 
$j_{\CV^\vee}^* \CF j_{\CV !} \iso \bar{\CR} (-1)$ in $\bar{D}(\Bbb P^\vee )$ (see 
 \cite{L} 1.6).
By \cite{L} 1.8,  one has a natural identification $$\CF i_{\CV !}  \iso \pi_{B_{\Bbb G_m}}^* [n], \tag 2$$ where $\pi_{B_{\Bbb G_m}} $ is the projection $ \CV^\vee \to B_{\Bbb G_m}$. Notice that $\pi_{B_{\Bbb G_m}*} \iso i^*_{\CV^\vee}$, so  $\pi_{B_{\Bbb G_m}}^*$ is left  adjoint to $ i^*_{\CV^\vee}$. Passing in (2) to the right adjoint functors, we get $$ i^!_{\CV }[n]\iso i^*_{\CV^\vee}\CF  . \tag 3$$

\enspace 

{\it Remark.} Other settings for Fourier transform of constructible sheaves 
  can be also used towards our aim (these are
 monodromic Fourier transform that identifies the subcategories of complexes with monodromic cohomology in $D(V)$ and $D(V^\vee )$,\footnote{Monodromic Fourier transform is the functor  $N\mapsto \text{holim}_a \, pr^\vee_! (pr^* N \otimes \mu^* j_{\Bbb A^1 *}\CL_a )[n+1]$, where $\ldots\twoheadrightarrow\CL_2 \twoheadrightarrow\CL_1$ are local systems on $\Bbb A^1 \smallsetminus \{ 0\}$ with unipotent Jordan block monodromy, rk  $\CL_a =a$. For the analytic version, see \cite{B} \S 6. }  and, for $\CD$-modules  or for $\ell$-adic sheaves in finite characteristic,  the  full Fourier transform that identifies  $D(V)$ with $D(V^\vee )$).

\subhead 4. Primitive cycles\endsubhead
Let $M$ be a non-constant irreducible perverse sheaf on $\Bbb P$. By the  theorem in 2, $\bar{\CR}(M)$ is an irreducible object of $\bar{\CM}(\Bbb P^\vee )$; let $M^\vee$ be the corresponding non-constant irreducible perverse sheaf on $\Bbb P^\vee$.
Let $c\in H^2 (\Bbb P, \Bbb Q(1))$ be the class of a hyperplane section. We have the primitive
decomposition\footnote{For an arbitrary irreducible $F$ this was proven in  \cite{D} (via \cite{BK} or \cite{G}) and \cite{M}.} $$\mathop\oplus\limits_{j\ge \text{max}\{a/2,0\}}   H^{a-2j}(\Bbb P, M(-j))^{\text{prim}}    \iso H^a (\Bbb P,M), \tag 4$$ where $H^{-i}(\Bbb P, M)^{\text{prim}}$ $:=
\Ker (c^i : H^{-i }(\Bbb P,M)\to H^{i}(\Bbb P,M)(i))$, $i\ge 0$,  the $j$-component of $\iso$ is multiplication by $c^j$. Set $H^a (\Bbb P, M)^{\text{coprim}}:=
\Ker (c : H^{a }(\Bbb P,M)\to H^{a+2}(\Bbb P,M)(1))$, $a\ge 0$, which equals component $j=2a$ of (4). Ditto for $M^\vee$.

\proclaim{Theorem} One has  canonical identifications $$ H^{a } (\Bbb P, M)^{\text{coprim}}\iso H^{a+2-n} (\Bbb P^\vee , M^\vee )^{\text{prim}}.\tag 5$$\endproclaim

\demo{Proof} The intermediate extension functor $j_{\CV !*} :\CM (\Bbb P)\to \CM (\CV )$, $j_{\CV !*} (M):= \text{Im} ({}^p H^0 j_{\CV !} (M)\to {}^p H^0 j_{\CV *} (M))$ identifies the category of irreducible perverse sheaves on $\Bbb P$ with that of those irreducible perverse sheaves on $\CV$ which are not supported on $\CV \smallsetminus \Bbb P =\{ 0\}/\Bbb G_m$.
 Since $\CF$ sends sheaves supported on $\CV \smallsetminus \Bbb P$ to constant sheaves and ${}^p H^0 j_{\CV !}(M)  = j_{\CV !*}(M)$, we see that (1) yields $j_{\CV^\vee}^* \CF j_{\CV !*} (M)\iso M^\vee (-1)$, hence  $j_{\CV^\vee !*}(M^\vee )= \CF j_{\CV !*} (M)(1)$. Applying (3), we get $i^!_{\CV}j_{\CV !*} (M)(1)=i^*_{\CV^\vee} j_{\CV^\vee !*}(M^\vee)[-n]$.
 Pulling it back by the smooth
projections  $\pi_{\CV}: V\to \CV$, $\pi_{\Bbb P}: V^{\circ}\to\Bbb P$ of relative dimension one, we get a canonical  isomorphism $$  i_{V}^! j_{V !*} (M^\flat )(1) \iso 
i_{V^\vee}^* j_{V^\vee !*} (M^{\vee \flat})[-n], \tag 6$$ where $M^\flat := \pi_{\Bbb P}^* M[1]$, $M^{\vee\flat} := \pi_{\Bbb P^\vee}^* M^\vee [1]$ are irreducible perverse sheaves on $V^\circ$, $V^{\vee \circ}$. Since $i_{V}^*$ is right t-exact and   $j_{V!*}(M^\flat )$ is irreducible, the complex $i^*_{V} j_{V!*}(M^\flat )$ is acyclic in degrees $\ge 0$; dually,  $i^!_{V} j_{V!*}(M^\flat )$ is acyclic in degrees $\le 0$. We get    (5) combining (6) with the next (well-known) lemma:

\proclaim{Lemma} There are canonical identifications $H^{a} i_{V}^* j_{V !*} (M^\flat ) \iso H^{a+1}(\Bbb P,M)^{\text{prim}}$, $  H^{a} i_{V}^! j_{V !*} (M) \iso H^{a-1}(\Bbb P,M(-1))^{\text{coprim}}.$
\endproclaim

\demo{Proof of Lemma} The canonical exact triangle $ i^!_V j_{V!*}(M^\flat )\to  i^*_V j_{V!*}(M^\flat )\to i^*_V j_{V*}(M^\flat )$ and the above acyclicity remark imply that
$$  i^!_V j_{V!*}(M^\flat )[1]=\tau_{ \ge 0} i^*_V j_{V*}(M^\flat ), \quad  i^*_V j_{V!*}(M^\flat )=\tau_{< 0}
i^*_V j_{V*}(M^\flat ). \tag 7$$ Now $i^*_V j_{V*}(M^\flat )\iso 
R\Gamma (V^{\circ} ,M^\flat )\iso R\Gamma (\Bbb P , M\otimes \pi_{\Bbb P *} \Bbb Q_{V\smallsetminus \{ 0\}} )[1]$, the first isomorphim comes since $M^\flat$ is $\Bbb G_m$-equivariant, the second is the projection formula.  Thus the evident exact triangle $\Bbb Q_{\Bbb P}\to 
\pi_{\Bbb P *} \Bbb Q_{V\smallsetminus \{ 0\}}\to \Bbb Q (-1)_{\Bbb P}[-1]$, its boundary map is $c$, yields isomorphism $ i^*_V j_{V*}(M^\flat )\iso \CC one ( R\Gamma (\Bbb P , M(-1))[-1]\buildrel{c}\over\to R\Gamma (\Bbb P , M)[1]). $ By (7) and (4), it provides the  identifications of the lemma, q.e.d. 
 \qed\enddemo\enddemo
  
\remark{Remark}  Only the case $a=0$ is needed for the aims of \cite{BFNP}. \endremark

\subhead 5. A description of $N_1$ \rm{(\cite{GG}, \cite{BFNP})}\endsubhead Let $X$ be an irreducible projective variety, $F$ be an irreducible perverse sheaf on $X$ whose support equals $X$ (the case we need is $F=\Bbb Q_X [\dim X ]$), $\CL$ be a very ample sheaf on $X$.  Let us describe the subspace $ N_1 H^0 (X, F)^{\text{prim}}$ of $H^0 (X, F)^{\text{prim}}$ (which is the intersection of kernels of all restriction maps to $H^0 (Y, F|_Y )$,  $Y$ is a closed proper subspace of $X$). 

We have the embedding $i_\CL : X\hra \Bbb P=\Bbb P_{\CL}$, $n= \dim H^0 (X,\CL)$, that corresponds to $\CL$; we assume that $X\neq \Bbb P$, so $M:=i_{\CL *}F$ is non-constant. Consider identification $\alpha : H^0 (X,F)^{\text{prim}} \iso H^{2-n}(\Bbb P^\vee , M^\vee )$ defined as the composition $H^0 (\Bbb P, M)^{\text{prim}}
= H^0 (\Bbb P, M)^{\text{coprim}}\iso H^{2-n}(\Bbb P^\vee , M^\vee )^{\text{prim}}
=H^{2-n}(\Bbb P^\vee , M^\vee )$ where $\iso$ is the isomorphism from (5).

For a constructible complex $G$ denote by $\CH^\cdot G$, $\tau_{\ge\cdot}$ its  usual ({\it not} perverse) cohomology sheaves and the canonical truncation. The projection $M^\vee \to \tau_{\ge 2-n}M^\vee$ yields
the map $H^{2-n}(\Bbb P^\vee , M^\vee )\to H^0 (\Bbb P^\vee , \CH^{2-n} M^\vee )$. Let $K_{\CL} \subset H^0 (X,F)^{\text{prim}}$ be the $\alpha$-preimage of its kernel. It coincides with the kernel of the composition $H^0 (X,F)^{\text{prim}}\hra H^0 (X,F)\buildrel{p_{(T)}^*}\over\lra H^{2-n}(\Bbb P^\vee ,\CR (M)) \to H^0 (\Bbb P^\vee , \CH^{2-n}\CR (M ))$,  which assigns to a primitive cycle the display of its images in $H^0 (Y, F|_Y )$ for 
 all hyperplane sections $Y$. (Indeed, by the decomposition theorem, $\CR (M)$ is the direct sum of $M^\vee$ and constant sheaves, so
   the kernel of a projection $\CH^{2-n}\CR (M )\to \CH^{2-n} (M^\vee )$ is a constant sheaf, and primitive cycles restrict to 0 on a general hyperplane section.)  Therefore $K_\CL$ consists of all primitive cycles whose restriction to each hyperplane section is 0.

Clearly $K_\CL \subset K_{\CL^{\otimes 2}}\subset $... Since every closed subscheme $Y\subset X$, $Y\neq X$, lies on a hypersurface of sufficiently high degree, we see that $N_1 H^0 (X, F)^{\text{prim}}$ equals $K_{\CL^{\otimes n}}$ for $n\gg 0$.

\end